\documentstyle[12pt]{article}

\begin{document}

\def\ii{\'{\char'20}}
\def\r{\rightarrow}
\def\err{\end{array}}
\def\bea{\begin{eqnarray}}
\def\eea{\end{eqnarray}}
\newcommand{\beq}{\begin{equation}}
\newcommand{\eeq}{\end{equation}}
\newcommand{\nn}{\nonumber}

\newtheorem{defn}{Definition}[section]
\newtheorem{lemma}[defn]{Lemma}
\newtheorem{theor}[defn]{Theorem}
\newtheorem{prop}[defn]{Proposition}


\begin{titlepage}

\title{\Large \bf A classification of fibre bundles over 
2-dimensional spaces\thanks{Talk given at the Meeting 
{\it New Developments in Algebraic Topology} 
(July 13-14, 1998, Faro, Portugal)}}

\author{Yu. A. Kubyshin\thanks{On leave of absence from the Institute
for Nuclear Physics, Moscow State University, 
119899 Moscow, Russia. E-mail address: kubyshin@theory.npi.msu.su}\\
\\
Department of Physics, U.C.E.H., University of the Algarve  \\
8000 Faro, Portugal}

\date{22 November 1999}

\maketitle

\begin{abstract}
The classification problem for principal fibre bundles over
two-dimensional CW-complexes is considered. Using the Postnikov
factorization for the base space of a universal bundle
a Puppe sequence that gives an implicit solution for the classification 
problem is constructed. In cases, when the structure group $G$
is path-connected or $\pi_{1}(G) = 0$, the
classification can be given in terms of cohomology groups.
\end{abstract}

\end{titlepage}


\section{Introduction}
\label{Sect.1}

In the present contribution we consider the classification problem for 
principal fibre bundles. We give a solution for the case
when $M$ is a two-dimensional path-connected CW-complex.

A motivation for this study came from calculations in two-dimensional
quantum Yang-Mills theories in Refs. \cite{AK1}, \cite{AK2}.
Consider a pure Yang-Mills (or pure
gauge) theory on a space-time manifold $M$ with gauge group $G$ 
which is usually assumed to be a compact semisimple Lie group. The vacuum
expectation value of, say, traced holonomy $T_{\gamma}(A)$ for a
closed path $\gamma$ in $M$ is given by the
following formal functional integral: 
\bea
   <T_{\gamma}> & = & \frac{1}{Z(0)} Z(\gamma),  \label{T-int} \\
   Z(\gamma) & = & \int_{\cal A} {\cal D}A \ e^{-{\cal S}(A)} \
   T_{\gamma}(A),          \label{Z-def}
\eea
where $A$ is a local 1-form on $M$, describing the gauge potential, and
${\cal S}(A)$ is the Yang-Mills action. The form $A$ is obtained
from the connection 1-form $w$ in a principal fibre bundle. Its 
base space is the space-time $M$ and its structure group is the gauge 
group $G$. We denote such bundles by  
$\xi = (E,M,G,p)$, where $E$ is called the total space and 
$p: E \rightarrow M$ is known as 
the projection map (see the definitions in Sect.~\ref{Sect2}).  
For a local cross-section $s$ of $\xi$ the 1-form $A$ is
obtained from $w$  as the pull-back $A=s^{*}w$ (see, for example, \cite{NS}, 
\cite{DFN}, \cite{D}).
Integration in Eq. (\ref{Z-def}) is performed over the space
${\cal A}$ of connections in $\xi = (E,M,G,p)$.
In quantum field theory often some heuristic
measure of integration is assumed, but in the two-dimensional case
the functional integral can be defined rigourously (see \cite{AL}). 
In general, the space
${\cal A}$ is not connected but may consist of a number of components
${\cal A}^{(\alpha)}$ labelled by elements $\alpha$ of some index set
${\cal B}$. Then the functional integral in (\ref{Z-def}) is given by a sum
over the elements of ${\cal B}$,
\beq
 Z(\gamma) = \sum_{\alpha \in {\cal B}} \int_{{\cal A}^{(\alpha)}}
 {\cal D}A^{(\alpha)} \ e^{-{\cal S}(A^{(\alpha)}} \
   T_{\gamma}(A^{(\alpha)}),          \label{Z-sum}
\eeq
each term being the functional integral over
connections in ${\cal A}^{(\alpha)}$. The set ${\cal B}$ 
of the components of ${\cal A}$ is in 1-1 correspondence with the 
set ${\cal B}_{G}(M)$ of non-equivalent
principal $G$-bundles over the manifold $M$. 

This feature has an analog in quantum mechanics. There functional
integrals are calculated over the space of paths which 
connect two given points, say $x_{0}$ and $x_{1}$, in a space-time
manifold $M$. Consider the case when $x_{0}=x_{1}$.
Then the space of integration is the space of loops based at
$x_{0} \in M$. We denote it by $\Omega M$. If the space-time
$M$ is multiply connected, then the space of loops splits into
components labelled by elements
\[
\alpha \in {\cal B} = \pi_{0}(\Omega M,*) \cong
\pi_{1}(M,x_{0}), 
\]
and the functional integral is given by a sum
of integrals over these components.

Similar to the case of quantum mechanics,
in gauge theory it is useful 
to have a characterization of the space ${\cal B}_{G}(M)$ in terms of
some objects that can be calculated relatively easy, e.g., in
terms of homotopy or cohomology groups. This is essentially the
classification problem. Usually two aspects
of this problem are considered:

1) Characterize the set ${\cal B}_{G}(M)$ of equivalence classes of 
principal fibre bundles with a given structure group $G$ over a 
given space $M$ in terms of some algebraic objects 
that are easy to calculate and that allow to
"enumerate" the bundles;

2) Given two non-equivalent bundles $\xi_{1}$ and $\xi_{2}$,
find a set of
characteristics that allows to distinguish them.

In the present article we consider the first aspect of the problem. 
We would like to mention that a powerful concept, giving a partial 
solution of the second aspect of the classification problem, is the 
concept of characteristic classes \cite{MilSta74}. 

With the motivation, coming from gauge theory, in the present article 
we address, in fact, a more general case, namely the 
classification of principal fibre bundles with structure group that    
is an arbitrary topological group and with base space that is 
a two-dimensional path-connected CW-complex. 
All necessary tools have been already developed 
in the literature. A method, which in many cases gives a solution 
of the classification problem and which we closely follow here, 
is discussed in the lectures by Avis and Isham \cite{AI}. 
Various particular cases were 
considered in previous works. For example, a classification
of principal fibre bundles with $G=U(1)$ over manifolds 
was studied in Refs. \cite{K}, \cite{AB}. 
In Ref. \cite{Witten2} a classification of principal fibre bundles 
over two-dimensional manifolds in the case when $G$ is a connected 
Lie group was obtained. There the classification was given in terms 
of elements of the group $\Gamma$, specifying the global
structure of $G$ through the relation $G = \tilde{G}/\Gamma$, where
$\tilde{G}$ is the universal covering group of $G$. 
However, we did not find in the literature a 
solution of the classification problem for $G$-bundles with more general 
base spaces $M$ and more general structure groups. In the present 
contribution we give a solution of this problem in terms of 
cohomology groups of $M$ in the case when $M$ is a 
two-dimensional CW-complex. Taking into account the character of the 
meeting, we tried to make our discussion rather pedagogical.  
To this end we gather together necessary definitions 
and explain in some detail main results and theorems from 
standard books on algebraic topology that are relevant for the 
classification problem and that we use in our derivation. 

The plan of the article is the following. In Sect.~\ref{Sect2} 
we describe the setting of the problem and recall main 
definitions. In Sect.~\ref{Sect3}
the Eilenberg-Maclane spaces and the Postnikov decomposition  
are discussed. In Sect.~\ref{Sect4} we construct a short exact sequence
which includes the set of equivalence classes of principal fibre 
bundles as an element.
This is the main result of the article.
Particular cases, when this set is characterized in terms
of a cohomology group, are also discussed. In the Appendix we 
list some basic definitions and examples of categories and functors and 
give a formulation of Brown's representation theorem. 

\section{Setting of the problem and main definitions}
\label{Sect2}

We begin with some definitions \cite{Span}, \cite{Swit}. 
Here we work in the category ${\cal PT}$ of pointed topological spaces
and base-point-preserving continuous maps (see the Appendix).  
Let $(E,e_{0})$ and $(X,x_{0})$ be pointed topological spaces. 
In order to simplify notations we omit the base points in some formulas 
if this does not cause confusion. 

\begin{defn}
\label{defn:hlp}
A map $p: E \rightarrow X$ is said to have a homotopy lifting 
property with respect to a topological space 
$Y$ if for every map $\tilde{f}: Y \rightarrow E$ and a homotopy 
$h: Y \times I \rightarrow X$ with $h(y,0) = p(\tilde{f}(y))$ 
there exists a homotopy $\tilde{h}: Y \times I \rightarrow E$
that covers $h$ (i.e., $p \circ \tilde{h} = h$) and  
such that $\tilde{h}(y,0)=\tilde{f}(y)$. 
\end{defn}

\begin{defn}
\label{defn:fibration}(see \cite{Span}, \cite{Swit}, \cite{Bred}, 
\cite{Whi})  
A map $p: E \rightarrow X$ is called a fibration if it has the homotopy
lifting property with respect to any topological space $Y$. 
A map $p: E \rightarrow X$ is called a weak fibration if it has the homotopy
lifting property with respect to all disks $D^{n}$, $n \geq 0$. 
\end{defn}

\noindent Space $F=p^{-1}(x_{0})$, where $x_{0}$ is the 
base point of $X$, is called the  fibre of the fibration $p$.
In general $p^{-1}(x)$ for $x \neq x_{0}$ need not be homeomorphic to 
$F$, but they all are of the same homotopy type.

Let $I=[0,1]$ be the closed unit interval. 

\begin{defn}
\label{defn:path}
The space of all continuous maps $w: \ I \rightarrow X$ with 
$w(0)=x_{0}$, topologized by the compact-open topology, is called the 
space of paths in $(X,x_{0})$ starting at $x_{0}$ and is denoted by  
$(PX,w_{0})$. 
\end{defn}

\noindent The base point of $(PX,w_{0})$ is the constant map given by 
$w_{0}(t)=x_{0}$ for all $t \in I$. 

\begin{defn}
\label{defn:loop}
The space of all continuous maps $w: \ I \rightarrow X$ with 
$w(0)=w(1)=x_{0}$, topologized by the compact-open topology, is called 
the space of loops in $(X,x_{0})$ based at $x_{0}$ and is denoted by 
$(\Omega X,w_{0})$.
\end{defn}

\noindent The base point $\omega_{0}$ is the trivial loop $w_{0}(t)=x_{0}$ 
for all $t \in I$. If $f$ is a map from $X$ to $Y$, then it induces 
the map $\Omega f: \Omega X \rightarrow \Omega Y$, defined 
in a natural way as $(\Omega f) (w) = f \circ w$ for $w \in \Omega X$. 
Higher loopings are defined  
by induction: $\Omega^{n} X = \Omega \left( \Omega^{n-1} X \right)$ and 
$\Omega^{n} f = \Omega \left( \Omega^{n-1} f \right)$. 

An important example of a fibration is the 
path-loop fibration of a pointed topological space $(X,x_{0})$.  
The fibration map $\tilde{p}: \ (PX,w_{0}) \rightarrow (X,x_{0})$ 
is defined by $\tilde{p}(w)=w(1)$ for a path $w$ in $X$. 

\begin{prop}
\label{prop:PX} \cite{Span}, \cite{Swit} 
The map $\tilde{p}: \ PX \rightarrow X$ is a fibration with 
fibre $\Omega X$. 
\end{prop}

\begin{defn}
\label{defn:fibbun}
A fibre bundle, denoted by $\xi = (E,B,F,p)$, is a collection consisting   
of a total space $(E,e_{0})$, a base space $(B,b_{0})$, a fibre 
$(F,e_{0})$ and a bundle map $p: (E,e_{0}) \rightarrow (B,b_{0})$ 
such that there exists an open covering $\{U\}$ of $B$ and for each 
open set $U \in \{ U \}$ there is a
homeomorphism $\phi_{U}: \  U \times F \rightarrow p^{-1}(U) \subset E$
such that the following property of local triviality holds:
\beq
 (p \circ \phi_{U}) (b,x) = p(\phi_{U}(b,x)) = b   \label{fb-triv}
\eeq
for all $(b,x) \in U \times F$. 
\end{defn}

\noindent A particular fibre over $b \in B$ is
$p^{-1}(b)$, and $p^{-1}(b)$ is homeomorphic to $F$ for all 
$b \in B$. We also have the inclusion map 
$i:\ (F,e_{0}) \hookrightarrow (E,e_{0})$.

Now we define the notion of a principal fibre bundle. 
Let $G$ be a topological group. 

\begin{defn}
\label{defn:prinfibbun} \cite{Swit}
A principal fibre bundle with structure group $G$ is a fibre bundle 
$\xi=(E,B,G,p)$ with fibre $G$ such that: 

(i) there exists an open covering $\{U\}$ of $B$, and for each open set 
$U \in \{ U \}$ there is a homeomorphism $\phi_{U}: U \times G 
\rightarrow p^{-1}(U)$, satisfying (\ref{fb-triv}); 

(ii) $G$ acts on $E$ with a right action:
$E \times G \rightarrow E$ in such a way that, for any $b \in U$ 
and $g \in G$, $\phi_{U}(b,g) = \phi_{U}(b,g_{0}) \cdot g$, where $g_{0}$ 
is the unit in $G$ and the dot denotes the right action of $G$ 
on $E$. 
\end{defn}

\noindent For a principal fibre bundle with structure group $G$ the term 
``principal $G$-bundle'' is often used. Note that under the action of $G$ 
points of $E$ move along the same fibre. Relations between fibre bundles 
and fibrations are given by the following propositions.   

\begin{prop}
\label{prop:fb-f2} \cite{Swit}
If $\xi=(E,B,F,p)$ is a fibre bundle, then $p: \  E \rightarrow B$ is a 
weak fibration. 
\end{prop}

\begin{prop}
\label{prop:fb-f1} \cite{Span}
Suppose that the base space of a fibre bundle is paracompact 
and Hausdorff; then such fibre bundle is a fibration. 
\end{prop}

It is clear that when we speak about the classification 
of principal fibre bundles we mean classification of classes 
of bundles defined by a natural equivalence relation. 

\begin{defn}
\label{defn:eqivfb}
Two fibre bundles $\xi_{1}=(E_{1},B,F,p_{1})$ and
$\xi_{2}=(E_{2},B,F,p_{2})$ with the same base and same fibre
are said to be equivalent 
(this is indicated by writing $\xi_{1} \cong \xi_{2}$) 
if there exists a
homeomorphism $h: E_{1} \rightarrow E_{2}$, satisfying the relation
$p_{2} \circ h = p_{1}$. If both bundles are principal fibre
bundles with structure group $G$, then $h$ has to satisfy 
the following additional property:
\[
  h(e \cdot g) = h(e) \cdot g, \; \; \; g \in G, \; \; e \in E.
\]
\end{defn}

\noindent The pointed set of all 
equivalence classes of principal $G$-bundles over $M$ is 
denoted by ${\cal B}_{G}(M)$. The base point of this set 
is the equivalence class of the trivial $G$-bundle 
$\xi_{0} = (M \times G, M, G, \mbox{pr}_{1})$ 
with the obvious right action of $G$ on $M \times G$. Here  
$\mbox{pr}_{1}: \ M \times G \rightarrow M$ is the projection onto the 
first factor, i.e., $\mbox{pr}_{1}(x,g)=x$ for $x \in M$ and $g \in G$. 

Now we recall the definition of an induced fibre bundle.
Consider a fibre bundle $\xi=(E,B,F,p)$ and a map
$f: M \rightarrow B$ from some space $M$. 

\begin{defn}
\label{defn:indfb}
The fibre bundle induced from $\xi$ by the map $f$, which is denoted by 
$f^{*}\xi$, is the fibre bundle with base $M$, fibre $F$,  
total space $E'$, given by    
\beq
E' = \left\{ (x,e) \in M \times E | f(x) = p(e) \right\} \subset 
M \times E, 
                  \label{ind-def}
\eeq
and bundle projection $p': \ E' \rightarrow M$ defined as 
$p'\left((x,e)\right) = x$. 
\end{defn}

\noindent In other words, the induced bundle $f^{*}\xi$ 
is constructed by glueing copies of the fibre $p^{-1}(b)$, $b \in B$, over
all  points $x$ of $M$ such that $f(x)=b$. There is also a map 
$f': \ E' \rightarrow E$ given by $f' \left( (x,e) \right) = e$.  
{}From the definition of the induced bundle it follows that 
\beq
    p \circ f' = f \circ p',    \label{ind-rel}
\eeq  
i.e., the following diagram is commutative: 

\vspace{0.2cm}
\begin{center}
\begin{tabular}{lcl}
   $E'$          & $\stackrel{f'}{\longrightarrow}$ & $E$  \\
 $\downarrow p'$ &     & $\downarrow p$  \\ 
   $M$           & $\stackrel{f}{\longrightarrow}$ & $B$ \\
\end{tabular}
\end{center}
\vspace{0.2cm} 

\noindent When $\xi$ is a principal $G$-bundle, the induced fibre bundle
$\xi' = f^{*}\xi$ is also a principal $G$-bundle with a $G$-action on
$E' \subset M \times E$ defined by $(x,e) \cdot g = (x,e \cdot g)$ 
for $x \in M$,  $e \in E$ and $g \in G$.

The following theorem is important for the 
classification of principal fibre bundles (see \cite{Swit}).

\begin{theor}
\label{theor:indfb}
Let $\xi=(E,B,G,p)$ be a principal $G$-bundle with base space being a 
CW-complex. Let $M$ be a CW-complex, and let $f_{1}, f_{2}$ be two maps 
from $M$ into $B$ that are homotopic, $f_{1} \simeq f_{2}$.
Then the induced principal $G$-bundles $f_{1}^{*}\xi$ and
$f_{2}^{*}\xi$ are equivalent, $f_{1}^{*}\xi \cong f_{2}^{*}\xi$.
\end{theor}

It turns out that for any topological group $G$ there exists a space 
$BG$, called its classifying space, and a principal $G$-bundle 
$\xi_{G} = (EG,BG,G,p_{G})$, called the universal $G$-bundle, 
such that every principal $G$-bundle $\xi=(E,M,G,p)$ is induced from 
$\xi_{G}$ by some map $f: M \rightarrow BG$ \cite{Steen}. Two homotopic maps 
$M \rightarrow BG$ induce equivalent bundles. The formal definition 
is the following. 

\begin{defn}
\label{defn:n-univ}
A bundle $\xi=(E,B,G,p)$ is $n$-universal if for any pointed CW-complex 
$(M,x_{0})$ with $\dim M = (n-1)$ there is a bijection 
$[M,x_{0};B,b_{0}] \cong {\cal B}_{G}(M)$.
A $G$-bundle is universal (or $\infty$-universal) provided it is 
$n$-universal for any $n \geq 1$; in this case the base space 
$B$ is called a classifying space of $G$. 
\end{defn}

\noindent The universal bundle and the classifying space are unique 
up to homotopy equivalence. 

The result quoted above can be derived in an elegant way if one considers 
${\cal B}_{G}(M)$ as a result of the action of a cofunctor ${\cal B}_{G}$ 
and uses Brown's representation theorem (see Theorem~\ref{theor:Brown} in 
the Appendix). We would like to discuss this issue in more detail.  

{}From now on we will be considering fibre bundles whose base spaces  
are path-connected pointed CW-complexes, but we do not restrict 
their dimension for the time being. 
Let $G$ be a topological group. Following Ref. \cite{Swit} let us 
introduce a cofunctor ${\cal B}_{G}: \  {\cal PW}' \rightarrow
{\cal PS}$ as follows (we recall that main definitions, some examples 
and the notation of categories and functors are listed in the Appendix; 
see in particular examples E.4 and E.8). For a pointed CW-complex
$(X,x_{0})$ its image ${\cal B}_{G}(X)$ is, as before, 
the pointed set of all equivalence classes of principal $G$-bundles 
over $X$. For a map $f: \  (X,x_{0}) \rightarrow (Y,y_{0})$, where 
$(X,x_{0})$ and $(Y,y_{0})$ are objects in ${\cal PW}'$, the morphism 
${\cal B}_{G}([f])$ is the map from ${\cal B}_{G}(Y)$ to 
${\cal B}_{G}(X)$ such that for a class $[\xi]$ of $G$-bundles over $Y$
the class ${\cal B}_{G}([f])[\xi] = [f^{*}\xi]$, where $f^{*}\xi$ is 
the induced $G$-bundle over $X$.

The important property of the cofunctor ${\cal B}_{G}$ is given by the
following theorem. 

\begin{theor}
\label{theor:BG} \cite{Swit} 
The cofunctor ${\cal B}_{G} : {\cal PW}' \rightarrow {\cal PS}$ satisfies 
the wedge axiom and the Mayer-Vietoris axiom. 
\end{theor}

\noindent These axioms are given in the Appendix (axioms W)
and  MV)). The property of ${\cal B}_{G}$, described by 
Theorem~\ref{theor:BG}, allows us to apply Brown's 
representation theorem (Theorem~\ref{theor:Brown} in the Appendix) 
and, as a consequence, obtain the result, 
we are interested in. Namely, Brown's theorem implies 
that for any $G$ there exists a
CW-complex $(BG,\tilde{b}_{0})$ and a principal $G$-bundle 
$\xi_{G}=(EG,BG,G,p_{G})$ such that for any pointed 
CW-complex $(M,x_{0})$ there is an equivalence 
\beq
  [M,x_{0};BG,\tilde{b}_{0}] \cong {\cal B}_{G}(M).   \label{MBG-B}
\eeq  
The correspondence is realized by the natural transformation $T_{\xi_{G}}$, 
constructed according to Example E.12 in the Appendix. In the case 
under consideration for any map 
$f:\ M \rightarrow BG \;$ we have $T_{\xi_{G}}(M)([f])=[f^{*}\xi_{G}]$. 
Brown's theorem states that $T_{\xi_{G}}$ is in fact a natural equivalence. 
In other words, homotopy classes of maps 
$f:\ M \rightarrow BG$ are in one-to-one correspondence with 
classes of equivalent bundles over $M$, so that a class $[f]$ 
corresponds to the class containing the induced
bundle $f^{*}\xi_{G}$. 
The space $(BG,\tilde{b}_{0})$ in (\ref{MBG-B}) is the classifying 
space for the cofunctor ${\cal B}_{G}$, while the bundle $\xi_{G}$ plays the 
role of the universal element of this cofunctor 
(see Definition~\ref{defn:univ}). The classifying space $BG$
and the universal element $\xi_{G}$ are unique up to homotopy 
equivalence. According to Definition~\ref{defn:n-univ} 
the $G$-bundle $\xi_{G}$ 
is the universal bundle. All $G$-bundles over $M$ can be 
obtained from it as induced bundles for various 
$f:M \rightarrow BG$.

An important property of the universal bundle is given 
by the following theorem. 

\begin{theor}
\label{theor:n-univ}  \cite{Steen} 
A principal $G$-bundle $\xi=(E,B,G,p)$ is $n$-universal 
if and only if it is path-connected and
$\pi_{q}(E)=0$ for $1 \leq q \leq n$.
A principal $G$-bundle $\xi=(E,B,G,p)$ is universal if and only if it is 
path-connected and $\pi_{q}(E)=0$ for all $q \geq 1$.
\end{theor}

Another property of $\xi_{G}$ can be obtained by considering the long
homotopy sequence of a bundle. We recall that for a 
fibre bundle $\xi = (E,B,F,p)$ the following exact sequence exists:
\beq
\ldots \stackrel{p_{*}}{\longrightarrow} \pi_{q+1}(B,b_{0})
\stackrel{\Delta}{\longrightarrow} \pi_{q}(F,e_{0})
\stackrel{i_{*}}{\longrightarrow} \pi_{q}(E,e_{0})
\stackrel{p_{*}}{\longrightarrow} \pi_{q}(B,b_{0})
\stackrel{\Delta}{\longrightarrow} \ldots,   \label{lhs}
\eeq
where the maps $p_{*}$ and $i_{*}$ are induced by the bundle projection
$p$ and the inclusion map $i: (F,e_{0}) \hookrightarrow (E,e_{0})$ 
respectively and $\Delta$ is constructed out of the boundary 
homomorphism in a standard way \cite{Span}, \cite{Swit}, \cite{Steen}.
Let us consider long homotopy sequence (\ref{lhs}) of  
the universal bundle $\xi_{G}$. Since all homotopy groups
of its total space $EG$ are trivial, the long homotopy sequence 
breaks into short exact sequences. The latter imply that 
\beq
    \pi_{q}(BG,\tilde{b}_{0}) \cong \pi_{q-1}(G,g_{0})   \label{BG-G}
\eeq
for $q \geq 1$ \cite{Steen}. We would like to note that the 
long homotopy sequence also exists for fibrations. 

As an immediate application of relation (\ref{BG-G}) we obtain the 
classification of principal $G$-bundles over the sphere $S^{n}$ 
in terms of the $(n-1)$th homotopy group of $G$:
\beq
{\cal B}_{G}(S^{n}) \cong [S^{n},x_{0};BG,\tilde{b}_{0}] \cong 
\pi_{n} (BG,\tilde{b}_{0}) \cong \pi_{n-1} (G,g_{0}).   \label{sphere}
\eeq

\section{Eilenberg-MacLane spaces and Postnikov decomposition}
\label{Sect3}

A general explicit construction of the universal $G$-bundle $\xi_{G}$ for any 
topological group $G$ is due to Milnor (see, for example, \cite{H}). 
For the classical groups $G=SO(k)$, $O(k)$, $SU(k)$, $U(k)$ and $Sp(k)$ 
there exist more convenient constructions of the universal bundles
with total spaces being Stiefel manifolds and classifying 
spaces being Grassmann manifolds \cite{DFN}, \cite{H}. 
However, in general it is not easy to describe the set of maps
$[M;BG]$. It turns out that in the case, when $M$ is of finite 
dimension, $BG$ can be substituted by some other space, which may 
be simpler to construct and to characterize. In order to 
discuss this it is necessary to first introduce 
the notion of $n$-equivalence.

\begin{defn}
\label{defn:n-equiv}
Let $(X,x_{0})$, $(Y,y_{0})$ be pointed topological spaces. 
A map $f: (X,x_{0}) \rightarrow (Y,y_{0})$ is called an $n$-equivalence 
($n \geq 1$) if it
induces a one-to-one correspondence between the path components of $X$
and $Y$ and for all $x_{0} \in X$ the induced map
\[
 f_{*}: \  \pi_{q}(X,x_{0}) \rightarrow \pi_{q}(Y,y_{0})
\]
is an isomorphism for $0 < q < n$ and an epimorphism for $q=n$.
A map $f: (X,x_{0}) \rightarrow (Y,y_{0})$ is called an 
$\infty$-equivalence or weak homotopy equivalence if it is an 
$n$-equivalence for all $n \geq 1$.
\end{defn}

\noindent For CW-complexes a map is
a homotopy equivalence if and only if it is a weak homotopy 
equivalence (see \cite{Swit}).
The importance of the notion of $n$-equivalence is seen from the 
following theorem. 

\begin{theor}
\label{theor:n-equiv} \cite{Swit} 
Let $(Y,y_{0})$ and $(Z,z_{0})$ be pointed topological spaces and let 
$f: (Y,y_{0})$ $\rightarrow (Z,z_{0})$ be an $n$-equivalence. 
Then for any pointed CW-complex $(M,x_{0})$ the induced map 
$f_{*}: [M,x_{0};Y,y_{0}] \rightarrow [M,x_{0};Z,z_{0}]$ 
is a bijection if $\dim M < n$ and is a surjection if $\dim M
\leq n$. 
\end{theor}

Applying this result to the classification of principal $G$-bundles over
a pointed CW-complex $(M,x_{0})$ with $\dim M \leq n$ we
conclude that if there exists a CW-complex $BG_{n}$ and a map
$f: BG \rightarrow BG_{n}$ which is an $(n+1)$-equivalence, then
\beq
 {\cal B}_{G}(M) \cong [M;BG] \cong [M; BG_{n}],    \label{BG-BGn1}
\eeq
i.e., all such bundles are in 1-1 correspondence with  
homotopy classes of maps $M \rightarrow BG_{n}$.

For the construction of $BG_{n}$ Eilenberg-MacLane spaces play an
important role. 

\begin{defn}
\label{defn:EML}
A space $X$, satisfying the properties: 
\bea
 \mbox{i)} & & \; \; X \; \; \; \mbox{is path connected};
 \nonumber \\
 \mbox{ii)} & & \pi_{q}\left( X,* \right) = \left\{
      \begin{array}{ll}
                 \pi, & \mbox{if} \; \;  q=n, \\
                  0 , & \mbox{if} \; \;  q \neq n, 
      \end{array} \right.    \nonumber
\eea
where $\pi$ is a group and $n$ is a positive integer, is called 
an Eilenberg-MacLane space and is denoted by $K(\pi,n)$. 
\end{defn}
 
For $n \geq 2$ the group $\pi$ in this definition 
must be abelian. It can be shown that for
any integer $n \geq 1$ and for any group $\pi$ (abelian if $n \geq 2$) 
the space $K(\pi,n)$ can always be constructed as a
CW-complex and is unique up to weak homotopy equivalence \cite{Span}, 
\cite{Swit}, \cite{Whi}.
Here follow some examples of Eilenberg-MacLane spaces.

\noindent (1) $K(Z,1) = S^{1}$.

\noindent (2) $K(Z,2) = CP^{\infty}$, where $CP^{\infty}$ is defined as
the direct limit $CP^{\infty}= \lim_{\rightarrow} CP^{n}$, i.e. the
CW-complex which is the union of the sequence $CP^{1} \subset CP^{2}
\subset \ldots$ of complex projective spaces,
topologized by the topology coherent with the collection
$\{ CP^{j} \}_{j \geq 1}$.

\noindent (3) $K(Z_{m},1) = L^{\infty}(m)$, the $\infty$-dimensional
lens space. In particular,
\beq
K(Z_{2},1) = RP^{\infty} = \lim_{\rightarrow} RP^{n}.
\label{KZ2}
\eeq

\vspace{0.3cm}

\noindent In general, Eilenberg-MacLane spaces are infinite 
dimensional, case 1 above is an exception. 

Let us discuss some properties of Eilenberg-MacLane
spaces which are needed for the sequel. We first recall the notion of
suspension. 

\begin{defn}
\label{defn:suspension}
A (reduced) suspension $(SX,z_{0})$ of a pointed topological space 
$(X,x_{0})$ is defined as the quotient
space of $X \times I$, where points of 
\[
   (X \times \{0\}) \cup (\{x_{0}\} \times I) \cup (X \times \{1\}) 
\] 
are indentified to a single point. 
\end{defn}

\noindent The base point $z_{0}$ is the image of
$(x_{0},0) \in X \times I$ under the quotient map $X \times I
\rightarrow SX$. It can be shown that if $X=S^{n}$, the $n$-dimensional 
sphere, then $S(S^{n})$ is homeomorphic to $S^{n+1}$ for all 
$n \geq 0$. 

It turns out that for any two pointed topological spaces 
$(X,x_{0})$, $(Y,y_{0})$ the sets $[Y,y_{0};\Omega X, w_{0}]$ and 
$[SX,z_{0};Y,y_{0}]$ possess group structures. 
This property follows from the fact that the loop space 
$(\Omega X,\omega_{0})$ is an $H$-group and the the suspension 
$(SX,z_{0})$ is an $H$-cogroup. In brief an $H$-group is a pointed
topological space with a binary operation, called a multiplication, 
which satisfies the group axioms up to homotopy. If a multiplication 
is homotopy commutative, then the $H$-group is said to be an abelian 
$H$-group. An $H$-cogroup is a pointed topological space with a 
co-multiplication that is homotopy associative, with a homotopy 
inverse, and with a homotopy identity map. An $H$-cogroup 
is called abelian if the co-multiplication is homotopy 
commutative. For 
rigorous definitions of $H$-groups and $H$-cogroups we 
refer the reader to the literature (see, for example, \cite{Span}, 
\cite{Swit}). The group structure on $[Y,y_{0};\Omega X,w_{0}]$ (resp. 
$[SX,z_{0};Y,y_{0}]$) is induced by the $H$-group structure of 
$(\Omega X,w_{0})$ (resp. $H$-cogroup structure of $(SX,z_{0})$). 

The following two propositions are relevant for us (see \cite{Span},
\cite{Swit}). 
Let $(X,x_{0})$ and $(Y,y_{0})$ be pointed topological spaces. 
\begin{prop}
\label{prop:adjcorr}
There is an isomorphism of groups
\[
  A: \ [SX,z_{0};Y,y_{0}] \rightarrow [X,x_{0};\Omega Y,\omega_{0}].
\]
\end{prop}
This isomorphism is called an adjoint correspondence. Let us 
remark that $\Omega$ can be viewed as a functor 
from the category ${\cal PT}$ of pointed topological spaces to the 
category of $H$-groups and continuous homomorphisms, whereas $S$ 
can be considered as a functor from ${\cal PT}$ to the 
category of $H$-cogroups and continuous homomorphisms. The 
existence of the isomorphism $A$ means that the 
functors $S$ and $\Omega$ are adjoint. 

\begin{prop}
\label{prop:SXOY}
The set $[SX,z_{0};\Omega Y, w_{0}]$ is an abelian group. 
\end{prop}
Using this it can be shown that for $n \geq 2$ 
the space $\Omega^{n}Y$ is an abelian 
$H$-group and $S^{n}X$ is an abelian $H$-cogroup. 
An immediate implication of Proposition~\ref{prop:SXOY} is that 
\[
  \pi_{n}(X,x_{0}) = [S^{n};X] \cong [S^{n-1};\Omega X] =
  \pi_{n-1}(\Omega X,w_{0}).
\]
Applying this result to the Eilenberg-MacLane space $K(\pi,n)$ 
($n \geq 1$) we see that
\[
  \pi_{q}(\Omega K(\pi,n),*) = \pi_{q+1}(K(\pi,n),*) = \left\{
      \begin{array}{ll}
                 \pi, & \mbox{if} \; \;  q=n-1, \\
                  0 , & \mbox{if} \; \;  q \neq n-1.
      \end{array} \right.
\]
It can be proved that these relations imply the existence of a weak 
homotopy equivalence between $K(\pi,n-1)$ and $\Omega K(\pi,n)$
(see \cite{Swit}). Therefore 
\beq
\Omega K(\pi,n) \simeq K(\pi,n-1).    \label{OmegaK-K}
\eeq
Hence, the set $[M;K(\pi,n)]$ carries a 
natural group structure. Moreover, since 
$K(\pi,n) \simeq \Omega^{2} K(\pi,n+2)$, according to Proposition 
\ref{prop:SXOY} this group is abelian.  
It can be specified further. Consider the
$n$th cohomology cofunctor $H^{n}$ with coefficients in $\pi$ 
(see Example E.11 in the Appendix). It turns out
that when restricted to the category ${\cal PW}'$ its classifying
space (see Definition~\ref{defn:univ}) is $Y=K(\pi,n)$ \cite{Span}. 
This means that there is a natural
equivalence between the cofunctors $\pi^{Y}$ and $H^{n}$, i.e., for any
$(M,x_{0}) \in {\cal PW}'$
\beq
  [M,x_{0};K(\pi,n),*]  \cong H^{n}(M;\pi).    \label{MK-Hn}
\eeq

The next important element, which we need for the classification of 
principal fibre bundles, is the Postnikov decomposition (called also 
the Postnikov factorization).
First let us recall the notion of a simple space. 

\begin{defn}
\label{defn:simple} 
A pointed topological space $(X,x_{0})$   
is called simple if $\pi_{1}(X,x_{0})$ acts trivially on homotopy groups
$\pi_{n}(X,x_{0})$ for all $x_{0} \in X$ and all $n \geq 1$. 
\end{defn}

\noindent Since $\pi_{1}(X,x_{0})$ acts on itself by conjugation, then for 
a path-connected $X$ simplicity implies that $\pi_{1}(X,x_{0})$ is abelian. 

In order to formulate the theorem on the Postnikov decomposition  
we first define some class of fibrations and the notion of a 
Postnikov system. 

\begin{defn}
\label{defn:prinfib}
Let $(B,b_{0})$, $(B',b'_{0})$ be pointed topological spaces,  
$\theta: (B,b_{0}) \rightarrow (B',b'_{0})$ a 
base-point-preserving map, and $p_{\theta}: E_{\theta} 
\rightarrow B$ the fibration, induced by $\theta$ from the 
path-loop fibration $p': PB' \rightarrow B'$. 
The fibration $p_{\theta}$ is called the principal fibration 
induced by $\theta$. If $B'$ is an Eilenberg-MacLane space 
$K(\pi,n)$, where $n \geq 1$ and $\pi$ is abelian, then 
$p_{\theta}$ is called a principal fibration of type $(\pi,n)$. 
\end{defn}

\noindent We recall that the path-loop fibration was defined in 
Sect.~\ref{Sect2}). The fibre of $p_{\theta}$ is $\Omega B'$. 

\begin{defn}
\label{defn:Moore-Post}
Let $(X,x_{0})$ and $(Y,y_{0})$ be pointed topological spaces, and let 
$f: (X,x_{0}) \rightarrow (Y,y_{0})$ be a map. A Moore-Postnikov 
decomposition of $f$ is a sequence 
\beq
\ldots \longrightarrow Y_{n+1} \stackrel{p_{n}}{\longrightarrow} Y_{n} 
\longrightarrow \ldots
\longrightarrow Y_{2} \stackrel{p_{1}}{\longrightarrow} Y_{1} 
\stackrel{p_{0}}{\longrightarrow} Y_{0}          \label{Post-seq}
\eeq
of spaces $Y_{n}$ and maps $f_{n}: X \rightarrow Y_{n}$ with the
following properties: 

\begin{enumerate}
\item for $n=0$ the space $Y_{0}=Y$ and $f_{0}=f$; 

\item for $n \geq 1$ the map $p_{n}: Y_{n+1} \rightarrow Y_{n}$ is 
a principal fibration of type $(\pi_{n}(X,x_{0}),n+1)$, and 
$p_{0}: Y_{1} \rightarrow Y_{0}=Y$ is a fibration; 

\item for $n \geq 0$ the map $f_{n+1}: X \rightarrow Y_{n+1}$ is a 
covering map of $f_{n}: X \rightarrow Y_{n}$, i.e. $f_{n} = 
p_{n} \circ f_{n+1}$;

\item for $n \geq 1$ the map $f_{n}$ is an $n$-equivalence; 

\item if $n \geq 1$ then $\pi_{q}(Y_{n}) = 0$ for $q \geq n$.
\end{enumerate}
\end{defn}

\begin{defn}
\label{defn:Post}
Let $(X,x_{0})$ be a pointed topological space. A Postnikov system  
(or Postnikov decomposition) of $X$ is the Moore-Postnikov decomposition 
of the map $f: X \rightarrow Y$, where $Y$ is the set of path-components 
of $X$ with the factor topology and $f$ is the natural projection. 
\end{defn}

\noindent Thus, if $X$ is path-connected, then $Y = \{y_{0}\}$ is the 
one-point space and $f$ is the constant map $f: X \rightarrow y_{0}$. 

\begin{theor}
\label{theor:Post} (see \cite{Span}, \cite{Bred}, \cite{Whi}, 
\cite{MoTa}) 
If $X$ is a simple pointed path-connected topological space, then there 
exists a Postnikov system of $X$. 
\end{theor}

\noindent Let $\pi_{n}$ denote the $n$th homotopy group $\pi_{n}(X,x_{0})$ 
of $X$. 
For every $n \geq 1$ the fibration 
$p_{n}: \ Y_{n+1} \rightarrow Y_{n}$ is induced by a map 
$\theta_{n}: \ Y_{n} \rightarrow K(\pi_{n},n+1)$ from 
the path-loop fibration 
$PK(\pi_{n},n+1) \rightarrow K(\pi_{n},n+1)$. Its fibre is
$\Omega K(\pi_{n},n+1) \simeq K(\pi_{n},n)$, see Eq. (\ref{OmegaK-K}).
Due to relation (\ref{MK-Hn}) the map 
$\theta_{n}$ corresponds to a characteristic class 
$k^{n+1}(X) \in H^{n+1}(Y_{n},\pi_{n})$ called the Postnikov invariant. 
It is obtained by the transgression of the 
so called fundamental class of some auxilliary fibration, which appears in 
this construction. We refer the reader to Refs. \cite{Span}, \cite{MoTa} 
for details.     

It can be shown that the spaces $Y_{n}$ of the Postnikov 
system are unique up to homotopy equivalence \cite{MoTa}, 
\cite{Whi}. A method to construct the Postnikov system of a 
given simple pointed path-connected CW-complex $X$ is described 
in Ref. \cite{Whi}. It appears also that a space can be reconstructed 
from its Postnikov system \cite{Whi}. Namely, suppose we are given a sequence 
${\cal P}=\{ \pi_{n},Y_{n-1},p_{n-1} | n \geq 1\}$, also called a 
Postnikov system (or fibred Postnikov system in \cite{Whi}), with the 
properties: 1) $\pi_{n}$ is an abelian group; 2) $Y_{0}$ is contractible; 
3) $p_{n}: Y_{n+1} \rightarrow Y_{n}$ is a fibration with fibre $F_{n}$ 
being an Eilenberg-MacLane space $K(\pi_{n},n)$; 4) the injection 
$\pi_{n}(F_{n}) \rightarrow \pi_{n}(Y_{n+1})$ is an isomorphism. Then the  
inverse limit $Y = \lim_{\leftarrow} Y_{n}$ is well-behaved. In general 
$Y$ does not have the homotopy type of a CW-complex. Nevertheless, there 
exists a CW-approximation $X_{\cal P}$ of $Y$ whose Postnikov system is 
${\cal P}$. In other words, if ${\cal P}$ is the Postnikov system 
of a space $X$, then $X$ and $X_{\cal P}$ are of the same homotopy type.  

The Postnikov decomposition will be used for the 
classification of principal $G$-bundles in the next section. There
$X$ will play the role of $BG$ and the spaces $Y_{n}$ will be $BG_{n-1}$.

\section{Classification of $G$-bundles over 2-dimensional CW - complexes}
\label{Sect4}

In this section we assume that $(M,x_{0})$ is a pointed path-connected
CW-complex of dimension $\dim M = 2$.

We apply Theorem \ref{theor:Post} to the space $BG$, 
the base space of the
universal $G$-bundle. Since $\dim M = 2$, for the 
classification of principal $G$-bundles over $M$ it suffices to
find a 3-equivalence $BG \rightarrow BG_{2}$.

As it was explained above, the total space $EG$ of the universal 
$G$-bundle is $\infty$-connected, hence the base $BG$ is path-connected. 
For the conditions of the Postnikov factorization theorem to hold we assume
that $(BG,\tilde{b}_{0})$ is simple. This implies, in particular, 
that $\pi_{1}(BG)$ is abelian. Let us discuss briefly this point.
Recall that for a topological group $G$ we have 
$\pi_{0}(G,g_{0}) = G/G_{0}$, where
$G_{0}$ is the path-connected component of the group unity, and
$\pi_{0}(G,g_{0})$ acts on higher homotopy groups of $G_{0}$
\cite{Steen}. According to relation (\ref{BG-G}) we have 
\beq
    \pi_{1}(BG,\tilde{b}_{0}) \cong \pi_{0}(G,g_{0}).  \label{BG1-G0}
\eeq
Then the condition of simplicity of $BG$ implies that
$\pi_{0}(G,g_{0})$ is abelian. One example, when this condition is 
trivially fulfilled, is $G$ being a path-connected group.
In this case $\pi_{1}(BG,\tilde{b}_{0}) \cong \pi_{0}(G,g_{0}) = 0$ 
and $BG$ is simple.
Another example is when $G$ is discrete 
and abelian. In this case $\pi_{0}(G,g_{0}) \cong (G,g_{0})$ and 
$\pi_{q}(G,g_{0}) = 0$ for all $q \geq 1$. 

In order to guarantee that $BG$ is simple we assume that $G$ is such that
$\pi_{0}(G,g_{0})$ is abelian and 
$\pi_{1}(BG,\tilde{b}_{0}) \cong \pi_{0}(G,g_{0})$ 
acts trivially on $\pi_{n}(G,g_{0})$ for $n \geq 1$. In addition we assume 
that $\pi_{0}(G,g_{0})$ is discrete. The latter condition will be used  
in the derivation of a short exact sequence which includes ${\cal B}_{G}(M)$. 

In what follows, $\pi_{n}$ denotes the $n$th homotopy group 
$\pi_{n}(BG,\tilde{b}_{0})$ of $BG$. 
The first non-trivial level of the
Postnikov diagram for $BG$ starts with $Y_{2}$:

\vspace{0.5cm}

\begin{center}
\begin{tabular}{cclcccc}
  & & $Y_{3}$ &
  $\stackrel{\theta_{2}'}{\longrightarrow}$ &
  $PK(\pi_{2},3)$ & $\hookleftarrow$ &
  $\Omega K(\pi_{2},3) \simeq K(\pi_{2},2)$ \\
  & & $\downarrow \  p_{2}$  &  &
  $\downarrow \tilde{p}_{2}$ & & \\
  & & $Y_{2}$ &
  $\stackrel{\theta_{2}}{\longrightarrow}$ & $K(\pi_{2},3)$ &  &   \\
  & & $\downarrow \  p_{1}$  &  & & & \\
  $BG$ & $\longrightarrow$ & $Y_{0}$ & 
 $\stackrel{\theta_{1}}{\longrightarrow}$ & $K(\pi_{1},2)$ & & \\
\end{tabular}
\end{center}

\vspace{0.5cm}
\noindent Recall that here $Y_{0} = \{ y_{0} \}$. 
 
Auxiliary path-loop fibrations are 
\[
\tilde{p}_{1}: PK(\pi_{1},2) \rightarrow K(\pi_{1},2) \; \; 
\mbox{with fiber} \; \; \Omega K(\pi_{1},2) \simeq K(\pi_{1},1)
\]
and
\[
\tilde{p}_{2}: PK(\pi_{2},3) \rightarrow K(\pi_{2},3) \; \; 
\mbox{with fiber} \; \; \Omega K(\pi_{2},3) \simeq K(\pi_{2},2).
\] 
The principal fibrations $p_{1}$ and $p_{2}$ are induced from them 
by $\theta_{1}^{*}$ and $\theta_{2}^{*}$ respectively.  
Of course, $Y_{2} = \{y_{0}\} \times \Omega 
K(\pi_{1},2) \simeq K(\pi_{1},1)$. The 3-equivalence, we
are looking for, is $f_{3}: BG \rightarrow Y_{3}$, so that $Y_{3}$ 
plays the role of $BG_{2}$. Its existence is guaranteed by the 
Postnikov decomposition Theorem \ref{theor:Post}.  
$f_{3}$ is related to the Postnikov invariant 
$k^{3}(BG) \in H^{3}(Y_{2}; \pi_{2})$, which corresponds to 
the map $\theta_{2}: Y_{2} \rightarrow K(\pi_{2},3)$ through the relation
$[Y_{2};K(\pi_{2},3)] \cong H^{3}(Y_{2}; \pi_{2})$ (see Sect.~\ref{Sect3}). 
Therefore, the non-trivial fibration, which is
important for us, is $p_{2}: Y_{3}=BG_{2} \rightarrow Y_{2} \simeq 
K(\pi_{1},1)$ with fiber $K(\pi_{2},2)$.

As the next step we build an exact sequence associated with
this bundle. This is a well known construction \cite{Span},
\cite{Swit}, \cite{Whi} (see also \cite{AI}) which we remind here.
Let $p: (E,e_{0}) \rightarrow (B,b_{0})$ be a fibration with fiber $F$, 
and let $i: (F,e_{0}) \hookrightarrow (E,e_{0})$ be the inclusion 
of the fibre. Firstly, we define
\[
\bar{E}= \left\{ (e,w) \in E \times PB | w(0)=p(e) \right\}
\subset E \times PB,
\]
where $(PB,w_{0})$ is the space of paths in $(B,b_{0})$.
Secondly, we define a path lifting function 
$\lambda: \bar{E} \rightarrow PE$ for $p$ 
(in Ref. \cite{Whi} this map is called a connection) by the
following conditions:
\[
\left( \lambda (e,w) \right)(0) = e \; \; \; \mbox{and} \; \; \;
p \left( \lambda (e,w) \right)  = w.
\]
Thus, $\lambda (e,w)$ is a path in $E$ which starts at $e \in E$ and
covers the path $w$ in $B$. If $w$ is a closed path in $B$, then
$\left( \lambda (e,w)\right) (1)$ is in the same fiber as $e$.
In a certain sense a lifting function is a generalization 
of the notion of a connection in differential geometry. 

\begin{prop}
\label{prop:plf} \cite{Span}, \cite{Swit}, \cite{Whi} 
A map $p: E \rightarrow B$ is a fibration if and only if 
there is a path lifting function for $p$. 
\end{prop}
Let us define a map $\rho: \Omega B \rightarrow F$ by the 
following formula:  
\beq
\rho (w) = \left( \lambda (e_{0},w^{-1}) \right) (1), \; \;
w \in \Omega B.    \label{rho-def}
\eeq
The map $\rho$ allows us to define a left exact sequence, 
characterized by Theorem~\ref{theor:seqn}, for the fibration 
$p$. 

\begin{defn}
\label{defn:leftex}
A sequence 
\[
\ldots \longrightarrow (Y_{n},y_{n}) \stackrel{f_{n}}{\longrightarrow} 
(Y_{n-1},y_{n-1}) \longrightarrow \ldots 
\ldots \longrightarrow (Y_{1},y_{1}) \stackrel{f_{1}}{\longrightarrow} 
(Y_{0},y_{0}) 
\]
of pointed topological spaces is called left exact if for any 
pointed topological space $(X,x_{0})$ the sequence 
\bea  
 & \ldots & \longrightarrow [X,x_{0};Y_{n},y_{n}] 
\stackrel{(f_{n})_{*}}{\longrightarrow} 
[X,x_{0};Y_{n-1},y_{n-1}] \longrightarrow \ldots  \nonumber \\ 
 & \ldots & \longrightarrow [X,x_{0};Y_{1},y_{1}] 
\stackrel{(f_{1})_{*}}{\longrightarrow} 
[X,x_{0};Y_{0},y_{0}]  \nonumber  
\eea
is exact in the category of pointed sets. 
\end{defn}

\begin{theor}
\label{theor:seqn} \cite{Swit}, \cite{Whi}
If $p: (E,e_{0}) \rightarrow (B,b_{0})$ is a fibration with 
fibre $(F,e_{0})$, and $\rho: \Omega B \rightarrow F$ is the map 
defined by Eq. (\ref{rho-def}), then the sequence  
\bea
\ldots \longrightarrow & \Omega^{n+1}B & 
       \stackrel{\Omega^{n} \rho}{\longrightarrow} \Omega^{n} F
       \stackrel{\Omega^{n} i}{\longrightarrow} \Omega^{n} E
       \stackrel{\Omega^{n}p}{\longrightarrow} \Omega^{n} B \longrightarrow
                     \nonumber  \\
\ldots & \longrightarrow & \Omega B
    \stackrel{\rho}{\longrightarrow} F \stackrel{i}{\longrightarrow} E
    \stackrel{p}{\longrightarrow} B.       \label{Omega-seq}
\eea
is left exact. 
\end{theor}
Sequence (\ref{Omega-seq}) is called the Puppe sequence of 
the fibration $p: E \rightarrow B$ \cite{MoTa}, \cite{P}. 

Now let us consider the Puppe sequence of the principal fibration
$p_{2}: Y_{3} \rightarrow Y_{2}$ with fiber $F \simeq K(\pi_{2},2)$ and
apply $[M,x_{0};-]$ to it. As it follows from Theorem \ref{theor:seqn},  
the sequence of pointed sets 
\bea  
\ldots & \longrightarrow & [M,x_{0};\Omega^{n+1} Y_{2},*] 
\stackrel{(\Omega^{n}\rho)_{*}}{\longrightarrow} 
[M,x_{0};\Omega^{n} K(\pi_{2},2),*] 
\stackrel{(\Omega^{n}i)_{*}}{\longrightarrow} 
[M,x_{0};\Omega^{n} Y_{3},*] \nonumber \\
& \stackrel{(\Omega^{n} p_{2} )_{*}}{\longrightarrow} & 
[M,x_{0};\Omega^{n} Y_{2},*] \longrightarrow \ldots 
\ldots \longrightarrow  
[M,x_{0}; Y_{3},*] 
\stackrel{(p_{2})_{*}}{\longrightarrow}
[M,x_{0}; Y_{2},*]  \label{M-Omega-seqn}
\eea
is exact. Sequences of this type sometimes are also called Puppe 
sequences.  

We extend this sequence to the right by adding 
$(\theta_{2})_{*}: [M;Y_{2}] \rightarrow [M;K(\pi_{2},3)]$, induced by
$\theta_{2}$. The extended sequence is exact too. The proof of this 
is rather simple and follows from properties of principal fibrations 
(see, for example, Ref. \cite{Bred}, \cite{P}). Indeed, let 
$\tilde{p}: PK \rightarrow K$ be a path-loop fibration over 
a space $(K,k_{0})$, and consider a map $\theta: B \rightarrow K$ 
and the fibration $p: E_{\theta} \rightarrow B$ induced from 
$\tilde{p}$ by $\theta$. Let $M$ be a pointed topological space. 
Any map $f: M \rightarrow E_{\theta}$ gives rise 
to a pair of maps: 1) $f_{1}=p \circ f: M \rightarrow B$ and 2) 
$f_{2}=\theta' \circ f: M \rightarrow PK$, where 
$\theta': E_{\theta} \rightarrow PK$ is the canonical map appearing in 
the definition of the induced fibration (see (\ref{ind-rel}) for an 
analogous map in the case of the induced bundle). Of course, 
$\theta \circ f_{1} = \tilde{p} \circ f_{2}$. 
Let $p_{*}$ and $\theta_{*}$ be the maps 
$p_{*}: [M;E_{\theta}] \rightarrow [M;B]$ and 
$\theta_{*}: [M;B] \rightarrow [M;K]$ induced by $p$ and $\theta$  
respectively. For any $x \in M$ 
$f_{2}(x)$ is a path in $K$ connecting the base point $k_{0}$ 
with the point 
\[
\left(f_{2}(x)\right)(1) = \tilde{p}\left(f_{2}(x)\right) = 
\left( \theta \circ f_{1} \right)(x) \in K. 
\]
Hence, the map $\theta \circ f_{1} = \theta_{*} f_{1}$ is 
homotopically trivial. Since $f_{1} = p_{*} f$, we conclude that 
$\mbox{im} \; p_{*} \subset \ker \; \theta_{*}$. The inverse inclusion is also 
easy to prove. Consider a map $f_{1}: M \rightarrow B$ such that 
$\theta_{*} f_{1} = \theta \circ f_{1}$ is homotopically trivial. 
This means that there is a homotopy $h: M \times I \rightarrow K$ with 
$h(x,0)=k_{0}$ and $h(x,1) = (\theta_{*} f_{1})(x)$ for all $x \in M$. 
It defines a path $w_{x} = h(x,-)$ in $K$, connecting $w_{x}(0)=k_{0}$ with 
\beq
  w_{x}(1) = (\theta_{*} f_{1})(x),  \label{final}
\eeq
and allows to introduce a function $f_{2}: M \rightarrow PK$ by 
the formula $f_{2}(x) = w_{x}$. Next we define a map 
$f: M \rightarrow E_{\theta} \subset B \times PK$ by 
$f(x) = (f_{1}(x), f_{2}(x))$. Relation 
(\ref{final}) gives that $\tilde{p} \left(f_{2}(x) \right) = w_{x}(1) = 
\theta \left( f_{1}(x) \right)$, and this guarantees that the image 
$f(x) \in E_{\theta}$. Since $p_{*} f = f_{1}$, we conclude that 
$\ker \; \theta_{*} \subset \mbox{im} \; p_{*}$. This proves that 
$\ker \; \theta_{*} =  \mbox{im} \; p_{*}$. 

Thus, for the principal fibration $p_{2}: Y_{3} \rightarrow Y_{2}$
we have the following extended Puppe sequence which is exact: 
\bea
\ldots & \longrightarrow & [M,x_{0};\Omega BG_{2},*]
   \stackrel{\left( \Omega p_{2} \right)_{*}}{\longrightarrow}
   [M,x_{0};\Omega K(\pi_{1},1),*]
   \stackrel{\left( \rho \right)_{*}}{\longrightarrow}
   [M,x_{0};K(\pi_{2},2),*]       \nonumber \\
   & \stackrel{\left( i \right)_{*}}{\longrightarrow} & [M,x_{0};BG_{2},*]
   \stackrel{\left( p_{2} \right)_{*}}{\longrightarrow} 
   [M,x_{0}; K(\pi_{1},1),*]
   \stackrel{\left( \theta_{2} \right)_{*}}{\longrightarrow} 
   [M,x_{0};K(\pi_{2},3),*].  
                           \label{M-Omega-seq}
\eea
Here we took into account that $Y_{3} = BG_{2}$ and 
$Y_{2} \simeq K(\pi_{1},1)$. Let us specify elements of this sequence. 

Firstly, isomorphism (\ref{MK-Hn}) gives that
\[ 
[M,x_{0};K(\pi_{1},1),*] \cong H^{1}(M;\pi_{1}) \; \; \; \mbox{and}
\; \; \; [M,x_{0};K(\pi_{2},2),*] \cong H^{2}(M;\pi_{2}).
\] 
Since the CW-complex $M$ is two-dimensional, the last term in the 
sequence 
\[
[M,x_{0};K(\pi_{2},3),*] \cong H^{3}(M;\pi_{2})=0.
\] 
This follows from the Universal Coefficient Theorem, see \cite{Span}, 
\cite{Swit}. 

Secondly, we use the relation 
\beq
    \Omega K(\pi,1) \simeq \pi    \label{OK-pi}
\eeq
valid for a discrete group $\pi$. To derive it take the CW-complex 
$X = K(\pi,1)$ with $\pi$ discrete. Consider the map 
$f: (\Omega X, *) \rightarrow \pi_{1}(X,*)$ which assigns to a loop 
$w \in \Omega X$ its homotopy class $[w] \in \pi_{1}(X,*)$, and a 
map $g: \pi_{1}(X,*) \rightarrow (\Omega X,*)$ which determines a 
representative in each class $[w] \in \pi_{1}(X,*)$. By definition 
$f$ is the homotopy inverse with respect to $g$ and vice versa, 
i.e., the map $f$ is a homotopy equivalence. Therefore, 
$\Omega K(\pi,1) \simeq \pi_{1}\left( K(\pi,1),* \right) \cong \pi$.    

Since, according to our assumption, $\pi_{1} = \pi_{1}(BG,\tilde{b}_{0}) 
\cong \pi_{0}(G,g_{0})$ is discrete, using result (\ref{OK-pi}) we 
conclude that 
\[
[M,x_{0};\Omega K(\pi_{1},1),*] \cong [M,x_{0};\pi_{1},*] = 0.
\]

Finally, we recognize the term $[M;BG_{2}]$ in sequence 
(\ref{M-Omega-seq}) as the
set that serves for the purpose of classification of principal
$G$-bundles over $M$ and that we are looking for. 
Indeed, according to the Postnikov
decomposition theorem there exists a 3-equivalence $f_{3}:BG
\rightarrow Y_{3}=BG_{2}$, so that relation (\ref{BG-BGn1}) with $n=2$ 
is valid. We have  
\[
   {\cal B}_{G}(M) \cong [M;BG] \cong [M;BG_{2}].
\]

Taking into account the properties of the terms in sequence
(\ref{M-Omega-seq}) and recalling that 
$\pi_{n}=\pi_{n}(BG,\tilde{b}_{0}) \cong \pi_{n-1}(G,g_{0})$, 
we arrive at the following result. 

\begin{theor}
\label{theor:main}
Let $(M,x_{0})$ be a path-connected pointed CW-complex of $\dim M =2$. 
Let $G$ be a group such that $\pi_{0}(G,g_{0})$ is abelian and 
discrete and acts trivially on the higher homotopy groups $\pi_{n}(G,g_{0})$ 
for $n \geq 1$. Then for the set ${\cal B}_{G}(M)$ of equivalence classes 
of principal $G$-bundles over $M$ there exists the following 
short exact sequence of pointed sets:  
\beq
0 \rightarrow H^{2}(M,x_{0};\pi_{1}(G,g_{0}),*) 
\rightarrow {\cal B}_{G}(M) \rightarrow 
H^{1}(M,x_{0};\pi_{0}(G,g_{0}),*) \rightarrow 0.  \label{HBH} 
\eeq 
\end{theor}
The maps in this sequence can be read from the extended Puppe 
sequence (\ref{M-Omega-seq}). 

This is our result on the classification of principal $G$-bundles. Further
characterization of ${\cal B}_{G}(M)$ requires the knowledge of
additional details about $M$ or $G$. 

Now we consider two particular cases and present
a discussion of the result.

\vspace{0.5cm}

1. \underline{$\pi_{1}(G,g_{0})=0$.}
It follows from (\ref{HBH}) that
\beq
   {\cal B}_{G}(M) \cong H^{1}(M;\pi_{0}(G,g_{0})).  \label{B-H1}
\eeq
This case includes the class of discrete groups. If $G$ is discrete,
$\pi_{0}(G) \cong G$ and the formula above simplifies further:
\beq
   {\cal B}_{G}(M) \cong H^{1}(M;G).  \label{B-H1a}
\eeq

2. \underline{$G$ is path-connected.}
Then $\pi_{0}(G,g_{0})=0$ and
\beq
   {\cal B}_{G}(M) \cong H^{2}(M;\pi_{1}(G,g_{0})).  \label{B-H2}
\eeq

\vspace{0.5cm}

The results in particular cases 1 and 2 above can be obtained, 
in fact, from the following general theorem.
\begin{theor} 
\cite{Bred} 
Let $(X,x_{0})$ be a pointed CW-complex, $(Y,y_{0})$ 
an $(n-1)$-connected
pointed space, and  $H^{q+1}(X;\pi_{q}(Y))=0=H^{q}(X;\pi_{q}(Y))$
for all $q \geq n$. Then there exists a one-to-one correspondence
\beq
  [X,x_{0};Y,y_{0}] \cong H^{n}(X;\pi_{n}(Y)).    \label{XY-Hn}
\eeq
\end{theor}

\noindent Let us apply this theorem to $X=M$ with $\dim M = 2$. In the 
case, when $Y = BG$
is 0-connected with $\pi_{2}(BG)=\pi_{1}(G)=0$, we obtain result
(\ref{B-H1}). If $Y=BG$ is 1-connected, i.e., 
$\pi_{1}(BG)=\pi_{0}(G)=0$, then (\ref{XY-Hn}) gives (\ref{B-H2}). The
latter case is also a result of the Hopf-Whitney classification
theorem \cite{Whi}.

Relations (\ref{B-H1}) and (\ref{B-H2}) become more concrete if
further properties of the group $G$ are known.
For completeness of the discussion, we present a list of the
first homotopy groups $\pi_{1}(G)$ for some connected Lie groups:

\begin{enumerate}

 \item Simply connected, $G = SU(n)$, $Sp(n)$: $\pi_{1}(G) = 0$.

 \item $G = SO(n)$, $n = 3$ and $n \geq 5$: $\pi_{1}(G) = Z_{2}$.

 \item $G=U(n)$: $\pi_{1}(G) = Z$.

\end{enumerate}

We would like to mention that for $G=U(1)$ actually a stronger result is
known \cite{AI}. Indeed, consider the Hopf bundle
$\xi_{n}=(S^{2n-1},CP^{n},U(1),p)$,
where the sphere is realized as
\[
  S^{2n-1} = \left\{ (z_{1}, \ldots ,z_{n}) \in C^{n} | \sum_{i=1}^{n}
  |z_{i}|^{2} = 1 \right\}
\]
and the bundle projection $p: S^{2n-1} \rightarrow CP^{n}$ is given by
\[
  p(z_{1}, \ldots ,z_{n}) = \left( \frac{z_{2}}{z_{1}}, \frac{z_{3}}{z_{1}},
  \ldots, \frac{z_{n}}{z_{1}} \right)
\]
for a neighbourhood corresponding to $z_{1} \neq 0$, etc. Since
$\pi_{q}(S^{2n-1}) = 0$ for $q < 2n-1$, the bundle is $(2n-2)$-universal.
Then one takes the direct limits $S^{\infty} = \lim_{\rightarrow} S^{n}$,
$CP^{\infty} = \lim_{\rightarrow} CP^{n}$. The bundle
$\xi_{\infty}=(S^{\infty},CP^{\infty},U(1),p)$ is
$\infty$-universal. Thus, $EU(1) = S^{\infty}$ and 
$BU(1) = CP^{\infty}$. This also can be obtained by using the 
Milnor construction \cite{H}. In Sect.~\ref{Sect3} we have already 
mentioned that $CP^{\infty} = K(Z,2)$. Then for a CW-complex of 
{\it any} dimension 
\beq 
{\cal B}_{U(1)}(M) \cong [M; BU(1)] \cong [M;CP^{\infty}] = 
[M; K(Z,2)] \cong H^{2}(M;Z).    \label{BG-U1}
\eeq
This is, of course, in agreement with our result (\ref{B-H2}).

Eq.~(\ref{B-H1a}), giving the classification of principal $G$-bundles with
discrete structure group, is in fact a particular case of a more
general result valid for CW-complexes $M$ of {\it any} dimension. 
It follows from the fact that for a discrete structure group $G$ 
\beq 
       BG \simeq K(G,1).     \label{K-Bpi}
\eeq
To show this, consider the universal covering
$UK(G,1) \rightarrow K(G,1)$ of $K(G,1)$. One can prove that it is
in fact a universal $G$-bundle \cite{Swit}, \cite{Steen}, \cite{P}. 
Indeed, $\pi_{1}(UK(G,1),*)=0$ by
definition, $\pi_{n}(UK(G,1),*)=\pi_{n}(K(G,1),*)=0$ for $n \geq 2$ 
and the group of automorphisms $Aut(UK(G,1)) \cong \pi_{1}(K(G,1),*)
\cong G$. Hence, $K(G,1)$ is the classifying space (see \cite{Steen} and 
Theorem \ref{theor:n-univ}), i.e., we have 
relation (\ref{K-Bpi}). With this we get 
\[
{\cal B}_{G}(M) \cong [M; BG] \cong [M; K(G,1)]
\cong H^{1}(M;G).                 \label{BG-disc}
\]
(cf. Eq.~(\ref{B-H1a})). In particular, the classifying space for 
$G=Z_{2}$ is $BZ_{2} \simeq K(Z_{2},1) = RP^{\infty}$, as one can see 
from (\ref{KZ2}). Using the Milnor construction it can be shown that 
the total space of the universal $Z_{2}$-bundle is $EZ_{2} = S^{\infty}$
\cite{H}. 

If the space $M$ is of certain type, 
formulas (\ref{HBH}), (\ref{B-H1}) and
(\ref{B-H2}) become more concrete. 
Thus, if $M$ is a two-dimensional differentiable manifold
one can use known formulas for cohomology groups. For example, if
$\pi$ is abelian then (see \cite{Bred})

\vspace{0.2cm}

\noindent 1) for $M$ compact and orientable 
\beq
H^{2}(M;\pi) \cong \pi,      \label{M-comp}  
\eeq

\noindent 2) for $M$ compact and not orientable 
\beq
H^{2}(M;\pi) \cong \pi/2\pi.     \label{M-noncomp}
\eeq
\vspace{0.2cm}

\noindent More expressions for cohomology groups of various 
two-dimensional surfaces can be found in \cite{GH}. 
In particular, for $M=S^{2}$ it is known that its first cohomology 
group $H^{1}(S^{2},\pi)=0$. Then from Theorem~\ref{theor:main}, 
Eq. (\ref{MK-Hn}) and Definition~\ref{defn:EML} it follows that 
\[
{\cal B}_{G}(M) \cong H^{2}(S^{2};\pi_{1}(G,g_{0})) \cong 
[S^{2};K(\pi_{1}(G),2)] \cong \pi_{2} \left( K(\pi_{1}(G),2),* \right) 
= \pi_{1}(G). 
\]
This is in accordance with Eq. (\ref{sphere}).  

If the structure group $G$ is path connected and simply connected,
i.e., $\pi_{0}(G,g_{0})=0=\pi_{1}(G,g_{0})$, 
then (\ref{HBH}) implies that the only (up to equivalence) 
existing principal $G$-bundle $\xi=(E,M,G,p)$ is the trivial one.

To finish our discussion of the result let us mention the classification
of principal fibre bundles over two-dimensional compact orientable
manifolds obtained by Witten in Ref. \cite{Witten2}.
It is known that any connected Lie group $G$ can be obtained 
as a quotient group $G = \tilde{G}/\Gamma$ (see, for example, \cite{BR}, 
\cite{FH}). 
Here $\tilde{G}$ is the unique (up to isomorphism) connected and simply 
connected Lie group, called the universal covering group. 
$\Gamma$ is a discrete subgroup of the center $Z(\tilde{G})$ of 
$\tilde{G}$. Witten showed that principal $G$-bundles over $M$ 
are classified by elements of $\Gamma$,
i.e. ${\cal B}_{G}(M) \cong \Gamma$. This agrees with Eq.~(\ref{B-H2}). 
Indeed, taking into account that 
$\pi_{1}(G,g_{0}) \cong \pi_{0}(\Gamma) \cong \Gamma$ and using 
Eq. (\ref{M-comp}), from (\ref{B-H2}) we get 
\[
{\cal B}_{G}(M) \cong H^{2}(M;\pi_{1}(G,g_{0})) 
\cong H^{2}(M;\Gamma) \cong \Gamma.
\]

\vspace{0.5cm}

\leftline{\bf Acknowledgements}

\vspace{0.25cm}

We would like to thank Manolo Asorey, Marco Mackaay, Jos\'e Mour\~ao
and Evgeny Troitsky for many fruitful discussions and useful remarks.
Financial support by Funda\c {c}\~ao para a
Ci\^encia e a Tecnologia (Portugal) under grants 
PRAXIS/2/2.1/FIS/286/ 94 and  
CERN/P/FIS/1203/98, from
fellowship PRAXIS XXI/BCC/4802/95 and 
from the Russian Foundation for Basic Research (grant
98-02-16769-a) are acknowledged.

\section*{Appendix}

\def\theequation{A.\arabic{equation}}
\setcounter{equation}{0}
\def\thedefn{A.\arabic{defn}}
\setcounter{defn}{0}
\def\thetheor{A.\arabic{theor}}

In this Appendix we give basic definitions and theorems, as well as 
some examples from the theory of categories, which we use in the article. 

\begin{defn}
\label{defn:cat}
A category ${\cal C}$ consists of 
\begin{enumerate}
\item a class of objects; 
\item sets of morphisms from $X$ to $Y$ for every ordered pair of objects 
$X$ and $Y$; such sets are denoted by $hom(X;Y)$; if $f \in hom(X;Y)$   
we write $f:X \rightarrow Y$; 
\item functions $hom(Y;Z) \times hom(X;Y) \rightarrow hom(X;Z)$ 
for every ordered triple of objects $(X,Y,Z)$; such a function is called a 
composition; if $f \in hom(X;Y)$ and $g \in hom(Y;Z)$, then the image 
of the pair $(g,f)$ in $hom(X;Z)$ is denoted by $g \circ f$. 
\end{enumerate}
Objects and morphisms satisfy the following two axioms: 

\begin{description}
\item[C1] Associativity: if $f \in hom(X;Y)$, $g \in hom(Y;Z)$ and 
$f \in hom(Z;W)$, then $h \circ (g \circ f) = (h \circ g) \circ f \in 
hom(X;W)$; 
\item[C2] Existence of an identity: 
for every object $X$ there exists a morphism 
$1_{X} \in hom(X;X)$ such that for every $g \in hom(Y;X)$ and every 
$h \in hom(X;Z)$, for all $Y,Z$, we have $1_{X} \circ g = g$ and 
$h \circ 1_{X} = h$.  
\end{description}  
\end{defn}

\noindent It can be shown that $1_{X}$ is unique. 
Here follow a few examples of categories. 

\begin{description}
\item[E.1] Category ${\cal S}$ of all sets and all functions.
\item[E.2] Category ${\cal T}$ of all topological spaces and all 
continuous maps.
\item[E.3] Category ${\cal G}$ of groups and homomorphisms.
\item[E.4] Category ${\cal PS}$ of pointed sets 
(sets with a distinguished element called a base point) 
and base-point-preserving functions. In this case 
the objects are denoted by $(X,x_{0})$, where $X$ is a 
set and $x_{0} \in X$ is 
the base point, and any $f \in hom(X,x_{0};Y,y_{0})$ is a 
function $f: X \rightarrow Y$ such that $f(x_{0}) = y_{0}$. 
\item[E.5] Category ${\cal PT}$ of pointed topological spaces and 
base-point-preserving continuous maps. 
\item[E.6] Category ${\cal T}'$ whose objects are topological spaces 
and $hom(X;Y) = [X;Y]$, the set of homotopy classes of continuous maps 
$f: X \rightarrow Y$. The homotopy class of a continuous map 
$f: X \rightarrow Y$ is denoted by $[f]$. For given $[f] \in hom(X;Y)$ and 
$[g] \in hom(Y;Z)$ the composition $[g] \circ [f]$ is defined as 
$[g \circ f]$.
\item[E.7] Category ${\cal PT}'$ whose objects are pointed topological 
spaces and $hom(X,x_{0};Y,y_{0}) = [X,x_{0};Y,y_{0}]$ is the set 
of homotopy classes of continuous base-point-preserving maps 
$f: (X,x_{0}) \rightarrow (Y,y_{0})$.  
\item[E.8] Homotopy category ${\cal PW}'$ whose objects are 
path-connected pointed CW-complexes and morphisms are homotopy classes of 
continuous base-point-preserving maps. Let us remind that a 
CW-complex is a topological 
space built up of cells, i.e., balls of various dimensions, glued together 
in a certain way. For a formal definition and detailed description of 
such spaces we refer the reader to Refs. \cite{Span}, \cite{Swit}, 
\cite{Whi}. 
\end{description}

\begin{defn}
\label{defn:equivcat}
Two objects $X$ and $Y$ in a category ${\cal C}$ are called equivalent 
if there exist morphisms $f \in hom(X;Y)$ and $g \in hom(Y;X)$ such 
that $g \circ f = 1_{X}$ and $f \circ g = 1_{Y}$. In this case the 
morphisms $f$ and $g$ are called equivalences. 
\end{defn}

\noindent We would like to mention that two sets are equivalent in 
the category ${\cal S}$ if and only if there exists a bijection between 
them. Two topological spaces $X$ and $Y$ are equivalent in ${\cal T}$ 
if and only if they are homeomorphic, and they are equivalent 
in ${\cal T}'$ if and only if they are homotopy equivalent. 

\begin{defn}
\label{defn:func}
A functor $F$ from a category ${\cal C}$ to a category ${\cal D}$ 
(we write $F: {\cal C} \rightarrow {\cal D}$) is a function that assigns   
\begin{enumerate}
\item to each object $X$ in ${\cal C}$ an 
object $F(X)$ in ${\cal D}$;
\item to each morphism $f \in hom_{\cal C}(X;Y)$ in ${\cal C}$ 
a morphism $F(f) \in hom_{\cal D}(F(X);F(Y))$  
such that the following axioms are fulfilled: 

\begin{description}
\item[F1] for each object $X$ in ${\cal C}$ $F\left(1_{X}\right) = 
1_{F(X)}$;
\item[F2] for every $f \in hom_{\cal C}(X;Y)$ and every 
$g \in hom_{\cal C}(Y;Z)$ in ${\cal C}$ $F(g \circ f) = F(g) \circ 
F(f) \in hom_{\cal D}(F(X);F(Z))$. 
\end{description}
\end{enumerate}
\end{defn}

\noindent The notion of a cofunctor is, in a certain sense, 
dual to that of a functor. 

\begin{defn}
\label{defn:cofunc}
A cofunctor $F^{*}$ from a category ${\cal C}$ to a category ${\cal D}$ 
(we write $F^{*}: {\cal C} \rightarrow {\cal D}$) is a function which 
assigns  
\begin{enumerate}
\item to each object $X$ in ${\cal C}$ an  
object $F^{*}(X)$ in ${\cal D}$;
\item to each morphism $f \in hom_{\cal C}(X;Y)$ in ${\cal C}$ 
a morphism $F^{*}(f) \in hom_{\cal D}(F^{*}(Y);F^{*}(X))$, 
such that the following axioms are fulfilled: 

\begin{description}
\item[CF1] for each object $X$ in ${\cal C}$ $F^{*}\left(1_{X}\right) = 
1_{F^{*}(X)}$;
\item[CF2] for every $f \in hom_{\cal C}(X;Y)$ and every 
$g \in hom_{\cal C}(Y;Z)$ in ${\cal C}$ $F^{*}(g \circ f) = F^{*}(f) \circ 
F^{*}(g) \in hom_{\cal D}(F^{*}(Z);F^{*}(X))$. 
\end{description}
\end{enumerate}
\end{defn}

Let us consider some examples. 

\begin{description}
\item[E.9] Let us fix a pointed topological space $(B,b_{0})$ in 
${\cal PT}$ and define the functor $\pi_{B}: {\cal PT} \rightarrow 
{\cal PS}$ as follows:
   \begin{description}
   \item[(i)] for each object $(X,x_{0})$ in ${\cal PT}$ the pointed set
       $\pi_{B}(X,x_{0})$ is  
       \[ 
               \pi_{B}(X,x_{0}) = [B,b_{0};X,x_{0}]; 
       \]
   \item[(ii)] for each morphism $f: (X,x_{0}) \rightarrow (Y,y_{0})$ in 
        ${\cal PT}$ the morphism $\pi_{B}(f)$ is given by
       \[
            \pi_{B}(f)[g] = [f \circ g] \in [B,b_{0};Y,y_{0}],
       \]
        where $[g] \in [B,b_{0};X,x_{0}]$. Often $\pi_{B}(f)$ is denoted 
        $f_{*}$ or $f_{\#}$. 
   \end{description}
The base point in the set $\pi_{B}(X,x_{0})$ is the homotopy class 
$[f_{x_{0}}]$ of the constant map $f_{x_{0}}: B \rightarrow x_{0} \in X$. 
\end{description}

\noindent Let $f$ and $f'$ be two continuous base-point-preserving 
maps $(X,x_{0}) \rightarrow (Y,y_{0})$ that are homotopic (written as  
$f \simeq f'$). It is easy to check that then $\pi_{B}(f) = \pi_{B}(f')$.
Hence $\pi_{B}$ is, in fact, a functor from the category ${\cal PT}'$ 
to ${\cal PS}$. 

\begin{description}
\item[E.10] For a fixed pointed topological space $(B,b_{0})$ in 
${\cal PT}$ we define the cofunctor $\pi^{B}: {\cal PT} \rightarrow 
{\cal PS}$ as follows:
   \begin{description}
   \item[(i)] for each object $(X,x_{0})$ in ${\cal PT}$ the pointed set 
       $\pi^{B}(X,x_{0})$ is  
       \[ 
               \pi^{B}(X,x_{0}) = [X,x_{0};B,b_{0}]; 
       \]
   \item[(ii)] for each morphism $f: (X,x_{0}) \rightarrow (Y,y_{0})$ in 
        ${\cal PT}$ the morphism $\pi^{B}(f)$ is given by
       \[
            \pi^{B}(f)[g] = [g \circ f] \in [X,x_{0};B,b_{0}],
       \]
        where $[g] \in [Y,y_{0};B,b_{0}]$. Often $\pi^{B}(f)$ is denoted 
        $f^{*}$ or $f^{\#}$. 
   \end{description}
The base point in $\pi^{B}(X,x_{0})$ is the homotopy class 
of the constant map $X \rightarrow b_{0} \in B$. 
\end{description}

\noindent For the same reason as in Example E.8 $\pi^{B}$ is in fact a 
cofunctor from the category ${\cal PT}'$ to ${\cal PS}$. 

\begin{description}
\item[E.11] Let $G$ be a group. The $n$th cohomology cofunctor 
$H^{n}$ with coefficients in $G$ is the cofunctor from ${\cal PT}'$ 
to ${\cal G}$ defined as follows: for any space $(X,x_{0})$ in ${\cal PT}'$ 
$H^{n}(X) = H^{n}(X,x_{0};G)$, the $n$th singular cohomology group with 
coefficients in $G$.
\end{description}  

\begin{defn}
\label{defn:nattran}
Let ${\cal C}$ and ${\cal D}$ be two categories and $F_{1}$ and 
$F_{2}$ be two functors from ${\cal C}$ to ${\cal D}$. A natural 
transformation $T$ from $F_{1}$ to $F_{2}$ is a function which 
assigns to each element $X$ in ${\cal C}$ a morphism 
$T(X) \in hom_{\cal D}(F_{1}(X);F_{2}(X))$ in ${\cal D}$ such that 
for each $f: X \rightarrow Y$ in ${\cal C}$ the 
following relation holds: 
\beq
  T(Y) \circ F_{1}(f) = F_{2}(f) \circ T(X).    \label{nt-def1} 
\eeq
\end{defn}

\noindent According to this relation the diagram 

\begin{center}
\begin{tabular}{ccc}
 $F_{1}(X)$ & $\stackrel{F_{1}(f)}{\longrightarrow}$ & $F_{1}(Y)$ \\
 $\downarrow \; T(X)$ &  & $\downarrow \; T(Y)$ \\
 $F_{2}(X)$ & $\stackrel{F_{2}(f)}{\longrightarrow}$ & $F_{2}(Y)$ \\
\end{tabular}
\end{center}

\noindent is commutative. The definition of the natural 
transformation between 
two cofunctors is obtained by inversion of the horizontal arrows 
in this diagram. Relation (\ref{nt-def1}) then takes the form 
\[
   T(X) \circ F^{*}_{1}(f) = F^{*}_{2}(f) \circ T(Y).  
\]

\begin{defn}
\label{defn:nateq}
Let ${\cal C}$ and ${\cal D}$ be two categories and let  
$F_{1},F_{2}: {\cal C} \rightarrow {\cal D}$ be two functors (cofunctors). 
A natural transformation $T$ from $F_{1}$ to $F_{2}$ is called 
a natural equivalence if $T(X)$ is an equivalence in the category 
${\cal D}$ for each object $X$ in ${\cal C}$. 
\end{defn}

\noindent The following example is important for us. 

\begin{description}
\item[E.12] Consider a cofunctor $H: {\cal PT}' \rightarrow {\cal PS}$ 
and fix a pointed topological space $(B,b_{0})$ and an element 
$u \in H(B)$. There is a standard way to introduce a natural 
transformation from $\pi^{B}$ to $H$. For every $(X,x_{0})$ in 
${\cal PT}'$ define a function $T_{u}(X): [X,x_{0};B,b_{0}] \rightarrow 
H(X)$ as follows: 
\[
T_{u}(X) \left([g] \right) = H \left([g]\right) (u) \in H(X)
\]
for any $[g] \in [X,x_{0};B,b_{0}]$. It is easy to check that $T_{u}$ 
is a natural transformation. 
\end{description}

The problem which arises in many applications is the following: for a 
given cofunctor $H$ find a space $(B,b_{0})$ and an element $u \in H(B)$ 
such that $T_{u}$ is a natural 
equivalence for the cofunctors $\pi^{B}$ and $H$. A solution is given 
by Brown's representation theorem discussed below. This theorem  
plays an important role for the classification of principal 
fibre bundles. 

An important class of cofunctors defined on ${\cal PW}'$ 
is formed by those cofunctors $H$ for which the following 
two axioms are fulfilled (see \cite{Swit}). 

\begin{description}
\item[W)] Wedge axiom. Let $\{ X_{\alpha} \}$ be an arbitrary family of 
$CW$-complexes in ${\cal PW}'$, and $\iota_{\alpha}: X_{\alpha} \rightarrow 
\vee_{\alpha} X_{\alpha}$ are the inclusions. 
Then for every $\alpha$ the morphism 
\[
H(\iota_{\alpha}): H\left( \vee_{\alpha} X_{\alpha}\right) 
\rightarrow \prod_{\alpha} H \left( X_{\alpha} \right) 
\]
is a bijection. 

\item[MV)] Mayer-Vietoris axiom. Consider a triple $(X;A_{1},A_{2})$, where 
$X$ is a CW-complex and $A_{1}$, $A_{2}$ are its subcomplexes such that 
$A_{1} \cup A_{2} = X$. Let $x_{1} \in H(A_{1})$ and $x_{2} \in H(A_{2})$ be 
arbitrary elements such that $H(i_{1})(x_{1}) = H(i_{2})(x_{2})$, where 
$i_{1}: A_{1} \cap A_{2} \rightarrow A_{1}$ and  
$i_{2}: A_{1} \cap A_{2} \rightarrow A_{2}$ are the inclusions. Then 
there exists an element $y \in H(X)$ such that $H(j_{1})(y) = x_{1}$ and 
$H(j_{2})(y) = x_{2}$, where $j_{1}: A_{1} \rightarrow X$ and 
$j_{2}: A_{2} \rightarrow X$ are the inclusions. 
\end{description}

\begin{prop}
\label{prop:piB} \cite{Swit}
For any $(B,b_{0})$ in ${\cal PW}'$ the cofunctor $\pi^{B} = [-;B,b_{0}]$ 
satisfies W) and MV).   
\end{prop}

To formulate Brown's representation theorem we need the notions  
of a universal element and of a classifying space. 

\begin{defn}
\label{defn:univ}
Let $H$ be a cofunctor on ${\cal PT}'$. An element $u \in H(B)$ for some 
$(B,b_{0})$ is called $n$-universal ($n \geq 1$) if 
\[
T_{u}(S^{q}): [S^{q},s_{0};B,b_{0}] \rightarrow H(S^{q})
\]
is an isomorphism for $1 \leq q < n$ and an epimorphism for $q = n$. 
An element $u \in H(B)$ is called universal if it is $n$-universal 
for all $n \geq 1$. In this case the space $B$ is called a 
classifying space for $H$. 
\end{defn}

\noindent 
For any cofunctor $H$ satisfying W) and MV) there exists a 
classifying space $(Y,y_{0})$ which is a CW-complex and a universal 
element $u \in H(Y)$. The classifying space for such a cofunctor is 
unique up to homotopy equivalence. Namely, if $(Y,y_{0}), (Y',y_{0}') \in 
{\cal PW}'$ are two classifying spaces for $H$ and $u \in H(Y)$, 
$u' \in H(Y')$ are corresponding universal elements, then there exists a 
homotopy equivalence $h: Y \rightarrow Y'$, unique up to homotopy, 
such that $H(h)(u')=u$ \cite{Span}, \cite{Swit}. 
Now we can formulate the theorem. 

\begin{theor}
\label{theor:Brown} (Brown's representation theorem, see \cite{Swit}) 
Let $H : {\cal PW}' \rightarrow {\cal PS}$ be a cofunctor satisfying 
axioms W) and MV). Then there exists a classifying space $(B,b_{0})$ in 
${\cal PW}'$ and a universal element $u \in H(B)$ such that 
\[
  T_{u}: \pi^{B} \rightarrow H
\]
is a natural equivalence. 
\end{theor}

\vspace{0.5cm}


\end{document}